\documentclass[11pt]{amsart}
\usepackage{graphicx}

\usepackage{amsmath,amsfonts,amsthm,amssymb,graphics,color, hyperref}

\newtheorem{theorem}{Theorem}
\newtheorem{prop}{Proposition}

\theoremstyle{remark}
\newtheorem{remark}{Remark}
\theoremstyle{definition}

\numberwithin{equation}{section}

\newcommand{\M}{\mathfrak{M}}

\newcommand{\pa}{\partial}
\newcommand{\eps}{\varepsilon}

\newcommand{\n}{\nabla}

\newcommand{\N}{\mathbb{N}}

\newcommand{\R}{\mathbb{R}}

\newcommand{\cL}{\mathcal{L}}

\renewcommand{\epsilon}{\varepsilon}

        \definecolor{pink}{rgb}{1,0,1}

\begin{document}
\title{Low eigenvalues and one-dimensional collapse} 

\author[Zhiqin Lu]{Zhiqin Lu}\address{Department of Mathematics, 410D Rowland Hall - University of California, Irvine, CA 92697-3875.} \email{zlu@uci.edu}

\author[Julie Rowlett]{Julie Rowlett} \address{Max Planck Institut f\"ur Mathematik, Vivatgasse 7, 53111 Bonn, Germany.} \email{rowlett@mpim-bonn.mpg.de}


\thanks{The first author is supported by NSF grant DMS-12-06748, and the second author acknowledges the support of the Max Planck Institut f\"ur Mathematik in Bonn and the Universit\"at G\"ottingen.}

\begin{abstract}
Our main result is that if a generic convex domain in $\R^n$ collapses to a domain in $\R^{n-1}$, then the difference between the first two Dirichlet eigenvalues of the Euclidean Laplacian, known as the fundamental gap, diverges.  The boundary of the domain need not be smooth, merely Lipschitz continuous.  To motivate the general case, we first prove the analogous result for triangular and polygonal domains.  In so doing, we prove that the  first two eigenvalues of triangular domains cannot be polyhomogeneous on the moduli space of triangles without blowing up a certain point.  Our results show that the gap generically diverges under one dimensional collapse and is bounded only if the domain is sufficiently close to a rectangle in two dimensions or a cylinder in higher dimensions.  
\end{abstract}
\maketitle

\tableofcontents

\section{Motivation and results} 
Our sign convention and boundary conditions for the Laplace equation on a domain $\Omega \subset \R^n$ shall be 
$$\Delta u =  \sum_{k=1} ^n \frac{\pa^2 u}{\pa x_k ^2} = - \lambda u, \quad \textrm{Dirichlet boundary condition:  }u = 0 \textrm{ on } \pa \Omega.$$
This equation arises from physics by separating the time and space variables in the wave equation.  Based on the physical interpretation of the eigenvalues, Marc Kac posed the now well known question, ``Can one hear the shape of a drum?'' \cite{kac}.  The mathematical formulation of this question is, if two Euclidean domains have the same set of Dirichlet eigenvalues, do they have the same shape?  Although the general answer to this question was proven to be ``One cannot hear the shape of a drum'' by Gordon, Webb, and Wolpert \cite{gww0, gww}, restricting to a specific moduli class of domains, it is possible for example to hear the shape of a triangle as proven by Durso \cite{durso}.  More recently, Grieser and Maronna have provided a simpler proof \cite{g-tri}.  Both proofs require the \em entire \em set of eigenvalues, but it is natural to speculate that the first three eigenvalues are sufficient to determine whether or not two triangles have the same shape.  Antunes and Freitas provided strong numerical evidence for this conjecture \cite{af-tri}.  

In some cases, a finite set of eigenvalues can indeed detect geometric features such as symmetry.  Polya proved that the first eigenvalue on a convex n-gon detects, up to scale, the regular n-gon for $n=3$ and $n=4$ \cite{polya1}.  For $n \geq 5$, the analogous result has not been proven and is known as Polya's Conjecture.  After the first eigenvalue the next natural object to study is the difference between the first eigenvalue and the rest of the spectrum, known as the \em fundamental gap. \em  The \em gap function \em is this difference multiplied with the square of the diameter of the domain, 
$$\xi(\Omega) := d^2(\Omega) \left( \lambda_2 (\Omega) - \lambda_1 (\Omega) \right),$$
for a bounded domain $\Omega \subset \R^n$.  If one scales a domain by a constant factor $c$, then the eigenvalues change according to : 
$$\lambda_k (c\Omega) = c^{-2} \lambda_k (\Omega).$$
Therefore the gap function is invariant under scaling of the domain and depends only on the shape of the domain.  

In \cite{tri-gap} we demonstrated that the gap function detects the equilateral triangle among all Euclidean triangles.  
\begin{theorem} [\cite{tri-gap}] Let $T$ be a Euclidean triangle.  Then 
$$\xi(T) \geq \frac{64 \pi^2}{9}$$
with equality if and only if $T$ is equilateral.  
\end{theorem} 

In a more significant paper, Andrews and Clutterbuck \cite{ac} proved that for any convex domain $\Omega \subset \R^n$, 
$$\xi(\Omega) \geq 3 \pi^2.$$

\subsection{Motivating Examples} 
There are very few domains for which one may explicitly compute the spectrum.  The first and most elementary example is a rectangle.  
\subsubsection{Rectangles} 
For a rectangle $R \cong [0, L] \times [0, W]$ the set of eigenvalues and eigenfunctions may be computed by separation of variables.  The eigenvalues are thus : 
$$\lambda_{j,k} = \frac{k^2 \pi^2}{L^2} + \frac{j^2 \pi^2}{W^2}.$$
Without loss of generality we assume $W \leq L$.  Then, the gap function 
$$\xi(R) = (L^2+W^2) \left( \frac{4 \pi^2}{L^2} + \frac{\pi^2}{W^2} - \frac{\pi^2}{L^2} - \frac{\pi^2}{W^2} \right) = \frac{3\pi^2 (L^2+W^2)}{L^2}.$$
We see immediately that the gap function detects the shape of a square, since for any rectangle $R$ 
$$\xi(R) \leq 6 \pi^2,$$
with equality if and only if $R$ is a square.  The gap function converges to the gap of a segment, $3 \pi^2$, if the rectangle collapses by letting $W \to 0$.  

\subsubsection{Circular Sectors} 
We next consider the Laplacian on a circular sector $S_\alpha$ of opening angle $\alpha \pi$ and radius $1$, which in polar coordinates $(r, \theta)$ is 
$$\Delta = \pa_r ^2 + r^{-1} \pa_r + r^{-2} \pa_\theta ^2, \quad (r, \theta) \in [0,1] \times [0, \alpha \pi].$$
Separating variables in the equation 
$$\Delta u = - \lambda u,$$
by assuming $u(r, \theta) = f(r) g(\theta)$, leads to the equations 
$$g''(\theta) = - \mu g(\theta), \quad g(0) = g(\alpha \pi) = 0, \quad \mu > 0,$$
$$r^2 f''(r) + rf'(r) + \lambda r^2 f(r) = \mu f(r).$$
The solutions to the first equation are 
$$g(\theta) = \sin(k \theta/\alpha), \quad k \in \N, \quad \mu = \mu_k = \frac{k^2}{\alpha^2}.$$
Re-writing the second equation as 
$$r^2 f''(r) + r f'(r) + (\lambda r^2 - \mu_k) f = 0, \quad f(0) = f(1) = 0,$$
one recognizes it as a Bessel equation, and the solutions are 
$$J_{k/\alpha} (\sqrt{\lambda} r), \quad J_{-k/\alpha} (\sqrt{\lambda} r),$$
where $J_x$ denotes the Bessel function of order $x$.  To satisfy the boundary conditions, the second solution is not allowed since it does not vanish at $r=0$, and thus 
$$\sqrt{\lambda} = \sqrt{\lambda_{j,k}} = Z_{j, k},$$
is the $j^{th}$ zero of the Bessel function of order $k/\alpha$.  The eigenvalues of the sector are therefore 
$$\lambda_{j,k} = Z_{j,k} ^2.$$ 
The above calculation is well known to spectral geometers but we have included it for the benefit of students and  researchers new to the field.  One of the aims of this paper is to accessibly present techniques and ideas which can be useful for eigenvalue problems and which require minimal, if any, sophisticated technical machinery.  

The following formulae are found in \cite{l-w, l-u}, and \cite{qw} for Bessel functions of real order and  \cite{watson} and \cite{olver} for Bessel functions of integer order.  The first and second zeros of the Bessel function of order $\nu$ are 
\begin{equation} \label{bzero} Z_{\nu, i} = \nu - \frac{a_i}{2^{1/3}} \nu^{1/3} + O(\nu^{-1/3}), \quad i = 1,2, \end{equation}
where $a_i$ is the $ i^{th}$ zero of the Airy function of the first kind so that   
\begin{equation} \label{airyzero}
a_1 \approx -2.33811, \textrm{ and } a_2 \approx -4.08795. \end{equation}
Consequently, we have the following estimates for the first two Dirichlet eigenvalues of the circular sector of opening angle $\alpha \pi$ and radius one
\begin{equation} \label{ev-sector} \lambda_i (S_\alpha) = \frac{1}{\alpha^2} + \frac{c_i}{\alpha^{4/3}} + O(\alpha^{-1}), \quad i=1,2,
\end{equation}
where
\begin{equation} \label{c} c_1 = -a_1 2^{2/3} \approx 3.71151827 \textrm{ and } c_2 = -a_2 2^{2/3} \approx 6.48921613. \end{equation} 
We therefore see that the gap function of a circular sector $S_\alpha$ of opening angle $\alpha \pi$ is asymptotic to 
$$\xi(S_\alpha ) = \frac{c_2 - c_1}{\alpha^{4/3}} + O(\alpha^{-1}), \quad \alpha \to 0.$$
In this case the gap function is unbounded as the sector collapses to the segment.  

These two simple examples show that the low eigenvalues, which determine the gap function, are sensitive to the geometry of convex planar domains which collapse to a segment.  In \cite{tri-gap}, we proved the following.    

\begin{theorem}[\cite{tri-gap}] \label{th:simplex}
Let $Y$ be an $n-1$ simplex for some $n \geq 2$.  Let $\{X_j \}_{j \in
\N}$ be a sequence of $n$ simplices each of which is a graph over $Y$.
 Assume the height of $X_j$ over $Y$ vanishes as $j \to \infty$.  Then
$\xi(X_j) \to \infty$ as $j \to \infty$.  More precisely, there is a
constant $C>0$ depending only on $n$ and $Y$ such that
$\xi(X_j)\geq C h(X_j)^{-4/3}$, where $h(X_j)$ is the height of $X_j$.
\end{theorem}

Since any triangle with unit diameter is a $2$-simplex with base 
$$Y = \{ (x,0) : 0 \leq x \leq 1\},$$
the above theorem shows that the gap function on a triangle which collapses to $Y$ blows up with the same asymptotic rate as the gap function on the sector $S_\alpha$ collapsing to $Y$.  In \S 2 below we explore different ways in which triangles may collapse to the segment and present an elementary proof independent of that in \cite{tri-gap} which shows that the gap function on any sequence of triangles which collapses to the unit interval is unbounded.  In \S 3 we prove that the gap function is sensitive to the geometry of collapsing polygonal  domains.  

We shall refer to the result we prove in \S 4 as a compactness theorem.  We show that for domains in $\R^{n}$ which collapse to a domain in $\R^{n-1}$, the gap function is sensitive to the geometry of the collapse.    Identify $\R^{n-1}$ with the subset
\[
\{(x_1,\cdots,x_{n-1},0)\mid (x_1,\cdots,x_{n-1})\in\R^{n-1}\} \subset \R^n.
\]

Let $\Omega\subset\mathbb R^n$ be a convex domain. Let
$E\subset\R^{n}$ be defined as
\[
E=\{(x_1,\cdots,x_{n-1},0)\mid \exists\, t, (x_1,\cdots, x_{n-1},t)\in\Omega\}.
\]

For $\eps>0$ define the one-parameter family of collapsing domains
\[
\Omega_\eps=\{(x_1,\cdots,x_n)\mid (x_1,\cdots, x_{n-1},\eps^{-1} x_n)\in\Omega\}.
\]

We call $E$ ``one-dimensional collapse'' and we denote
\[
E=\lim_{\eps\to 0} \Omega_\eps.
\]
Note that $E$ is canonically identified with a convex domain in $\mathbb R^{n-1}$.

For the sake of simplicity, we assume that $\Omega$ can be represented by
\[
\{(x,y)\mid x\in E, 0\leq y\leq h(x)\},
\]
where $h(x)$ is a  nonnegative Lipschitz continuous concave function over $E$.

\begin{theorem}\label{main} Let $\sigma=\max h(x)$.
Assume that 
\[
{\rm diam}\,\{x \in E \mid h(x)\geq\sigma-\delta\}\to 0
\]
as $\delta\to 0$.
Then
\[
\lambda_2(\Omega_\eps)-\lambda_1(\Omega_\eps)\to \infty
\]
as $\eps\to 0$.
\end{theorem}

A similar result to the above was proven by Borisov and Freitas \cite{bf-thin} who determined the asymptotics of both the eigenvalues and eigenfunctions for domains in $\R^d$ which collapse to $\R^{d-1}$.  The difference between our work and theirs is that they require the defining function, which corresponds to our function $h$ to be smooth in a neighborhood of its maximum point, whereas our function $h$ need only be Lipschitz continuous.  In particular, that work cannot be applied directly to simplices or triangles.  In dimension 2, Theorem \ref{th:simplex} follows from Friedlander and Solomyak's results which give a two-term asymptotic expansion for the eigenvalues in terms of a certain one-dimensional Schr\"odinger operator \cite{strip}.  We note however that the proofs presented here are significantly more elementary and independent to those arguments, and thus make a complementary contribution to the more refined and technical results and arguments in \cite{bf-thin} and \cite{strip}.  Eigenvalues of thin and collapsing domains are also relevant to physics; see \cite{gr1}, \cite{gr2}.  Further results include \cite{grj}.

\section{Collapsing triangular domains}
Recall the variational or mini-max principle for the first eigenvalue on a domain $\Omega \in \R^n$ 
\begin{equation} \label{r1} \lambda_1 = \inf_{f \in H^{1,2} _0 (\Omega)  } \frac{ \int_\Omega |\nabla f |^2 }{\int_\Omega f^2}.  \end{equation}
This is the infimum of the Rayleigh-Ritz quotient.  Above $H^{1,2} _0 (\Omega)$ is the set of functions $f \in \cL^2 (\Omega)$ such that 
$$\int_\Omega \left( |\nabla f|^2 + f^2 \right) < \infty, \quad f|_{\pa \Omega} =0$$
The second eigenvalue  is 
\begin{equation} \label{r2} \lambda_2 = \inf_{f \in H^{1,2} _0 (\Omega) : f \perp f_1} \frac{ \int_\Omega |\nabla f |^2 }{\int_\Omega f^2 dx dy}, \end{equation}
where $f_1$ is the eigenfunction for $\lambda_1,$ and orthogonality is with respect to $\cL^2.$  These formulae are in \cite{chavel}.  

An equivalent formula found in \cite{cour-hil} is the so-called ``maxi-min'' principle
\begin{equation} \label{maximin} \lambda_k = \inf_{L \subset H^{1,2} _0 (\Omega), \, \dim(L) = k} \left\{ \sup_{f \in L} \frac{ \int_\Omega |\nabla f |^2 }{\int_\Omega f^2 }  \right\}. \end{equation} 
The maxi-min principle can be used to prove domain monotonicity:  if $\Omega \subset \Omega'$ then 
$$\lambda_k (\Omega) \geq \lambda_k (\Omega').$$  

We shall use the following standard notation:  for functions $f, g: \R \to \R$ and $L \in \R \cup \infty$ we say that 
$$f(x) = O(g(x)) \textrm{ as } x \to L,$$
if there exists a constant  $C \in \R$ such that 
$$|f(x)| \leq Cg(x) \quad \textrm{ as } x \to L,$$
and we say that 
$$f(x) = o(g(x)) \textrm{ as } x \to L$$
if 
$$\lim_{x \to L} \frac{f(x)}{g(x)} = 0.$$
We write 
$$f(x) \sim g(x) \textrm{ as } x \to L$$
if there exist constants $c', c''$ such that 
$$c' g(x) \leq f(x) \leq c'' g(x), \quad \textrm{ as } x \to L,$$
and 
$$f(x) \approx g(x) \textrm{ as } x \to L$$
if 
$$\lim_{x \to L} \frac{f(x)}{g(x)} = 1.$$
In some arguments we will use the following constants\footnote{This constant arises from the asymptotic formula for the first two zeros of the derivative of the Bessel function which is related to the first two Dirichlet eigenvalue of an obtuse isosceles triangle; see \cite{freitas}.}  
\begin{equation} \label{c'} c_1' = -a_1 ' 2^{2/3} \approx 1.61722832, \quad c_2 ' = - a_2 ' 2^{2/3} \approx 5.15619226. \end{equation} 

\subsection{A representation of the moduli space of triangles}
\begin{figure}[ht]
\includegraphics[width=2in]{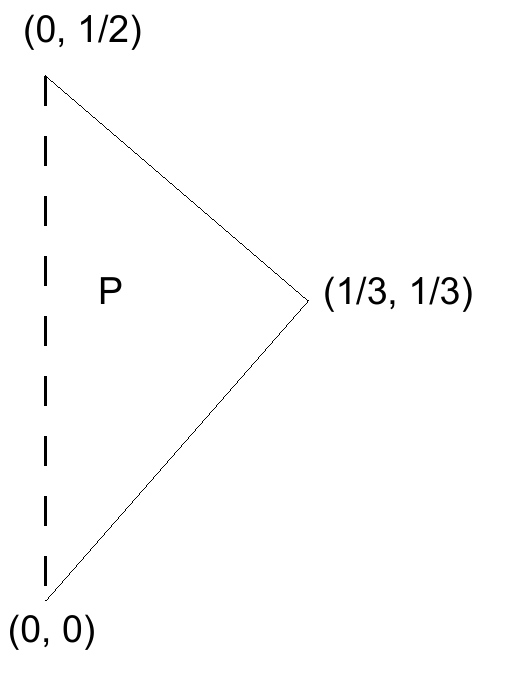}\caption{Moduli space of triangles.}\label{mod}
\end{figure}

Since the gap function is invariant under scaling, we restrict to triangles with diameter one.  Such a triangle has angles 
$$0 < \alpha \pi \leq \beta \pi \leq \pi - \alpha \pi - \beta \pi.$$
The moduli space of triangles, which we shall call $\M$, can thus be identified with a triangle in the $\alpha \times \beta$ plane; see Figure \ref{mod}.  We are interested in the behavior of $\xi$ approaching the boundary of $P$ which is the dashed vertical segment in Figure \ref{mod}; this corresponds to triangles which degenerate to a segment.  

\subsubsection{One collapsing angle}  
Consider the triangle $T$ with angles $0 < \alpha \pi \leq \beta \pi \leq \pi - \alpha \pi - \beta \pi,$ and assume for some fixed $\epsilon > 0,$ $\beta \geq \epsilon.$  Let the side opposite $\alpha \pi$ have length $A$, the side opposite $\beta \pi$ have length $B,$ and the third side have  length one.  The Law of Sines states that
$$\frac{\sin(\alpha \pi)}{A} = \frac{\sin(\beta \pi )}{B} = \frac{\sin(\pi - \alpha \pi - \beta \pi)}{1}.$$
Then 
$$B = \frac{\sin(\beta \pi )}{\sin(\alpha \pi + \beta \pi)} \approx 1 \textrm{ as } \alpha \to 0,$$
and 
$$B - 1 \sim \alpha \textrm{ as } \alpha \to 0.$$
For $\alpha$ small, we approximate the eigenvalues using two sectors, a larger sector $S$ of radius $1,$ and a smaller sector $\sigma$ of radius $B$ both with opening angle $\alpha \pi$.  By domain monotonicity,
$$\lambda_2 (T) \geq \lambda_2 (\textrm{large sector}) = \lambda_2 (S),$$
$$\lambda_1 (T) \leq \lambda_1 (\textrm{small sector}) = \lambda_1 (\sigma).$$
Then,
$$\lambda_2 (S) - \lambda_1 (\sigma) \leq \xi(T) \leq \lambda_2 (\sigma) - \lambda_1 (S).$$
Since 
$$\lambda_k (\sigma) = B^{-2} \lambda_k (S), \quad \forall k \in \N,$$
we have the estimate 
\begin{equation} \label{vp-sect} \lambda_2 (S) - B^{-2} \lambda_1 (S) \leq \xi (T) \leq B^{-2} \lambda_2 (S) - \lambda_1 (S).  \end{equation} 
Since $B -1 \sim \alpha$ and $\lambda_k (S) \approx \frac{1}{\alpha^2}$ as $\alpha \to 0$, there are constants $c, c' > 0$ such that 
$$\xi(S) - \frac{c}{\alpha} \leq \xi (T) \leq \xi(S) +\frac{c'}{\alpha},$$
and by (\ref{ev-sector}) 
$$\xi(S) = \frac{c_2 - c_1}{\alpha^{4/3}} \implies \lim_{\alpha \to 0} \frac{ \xi(T)}{\xi(S)} =1.$$


\subsubsection{Two collapsing angles:  different rates}
In fact, the same argument can be used to prove the following.  

\begin{prop} \label{prop-sector}  Let $\{T_n\}$ be a sequence of triangles with diameter one and angles $\alpha_n \pi$, $\beta_n \pi$ and $\pi(1 - \alpha_n - \beta_n)$.  Let $S_n$ be the sector with opening angle $\alpha_n \pi$ and radius $1$.  Assume that $\alpha_n < \beta_n$ and $\beta_n \to 0$ as $n \to \infty$.  Then 
$$\xi(T_n) \approx \xi (S_n) = \frac{c_2 - c_1}{\alpha_n ^{4/3}} \quad \textrm{ as } n \to \infty \iff \alpha_n = o(\beta_n ^3) \textrm{ as } \beta_n \to 0.$$
\end{prop} 
\begin{proof}  For simplicity, we shall abuse notation and drop the subscript $n$.  By the Law of Sines as above, 
$$B = \frac{\sin(\beta \pi)}{\sin(\alpha \pi + \beta \pi)}.$$
Estimating with two sectors $S$ and $\sigma$ of opening angle $\alpha \pi$ and radii $1$ and $B$, respectively, we again have the estimate (\ref{vp-sect}).  In this case since both $\alpha$ and $\beta$ tend to $0$, we estimate 
$$B^{-2} \approx 1, \quad B-1 \sim \frac{\alpha}{\beta}.$$
Therefore (\ref{vp-sect}) becomes 
$$\xi(S) - \frac{c \alpha}{\beta} \lambda_1 (S) \leq \xi (T) \leq \xi(S) + \frac{c' \alpha}{\beta} \lambda_2 (S),$$
for some positive constants $c$ and $c'$.  By (\ref{ev-sector}) for some positive constants $C$ and $C'$ we have 
$$\xi(S) - \frac{C}{\alpha \beta} \leq \xi (T) \leq \xi(S) + \frac{C'}{\alpha \beta}.$$
By (\ref{ev-sector}), 
$$\frac{\xi(T)}{\xi(S)} \to 1 \iff \frac{\alpha^{1/3}}{\beta} \to 0 \quad \textrm{ as } \beta \to 0,$$
which is equivalent to 
$$\alpha = o (\beta^3), \quad \textrm{as } \beta \to 0.$$
\end{proof} 

\begin{remark} We shall see below that if two angles collapse at sufficiently similar rates, then the gap function is asymptotic to 
$$\frac{4(c_2' - c_1')}{\alpha^{4/3}}, \quad \textrm{ as } \alpha \to 0.$$
An interesting open problem is to prove that the eigenvalues of triangles are polyhomogeneous functions on the moduli space of triangles.  This is not quite the full story, because as one can already see from our results, the gap function has different limits approaching the point $(\alpha, \beta) = (0, 0)$ along different trajectories.  Therefore, one must blow up the point $(\alpha, \beta) = (0,0)$, and attempt to prove polyhomogeneity on a suitably blown up space.  This is currently under investigation by Grieser and Melrose \cite{gm-phg}.  
\end{remark} 

\subsubsection{Two collapsing angles:  similar rates}
Consider obtuse isosceles triangles with diameter one.  By \cite{freitas}, 
$$\lambda_1 (T) = \frac{4}{\alpha^2} + \frac{4 c_1'}{\alpha^{4/3}} + O(\alpha^{-2/3}),$$
$$\lambda_2(T) = \frac{4}{\alpha^2} + \frac{4 c_2'}{\alpha^{4/3}} + O(\alpha^{-2/3}).$$
Then, 
$$\xi(T) \approx \frac{4(c_2'-c_1')}{ \alpha^{4/3}} + O(\alpha^{-2/3})\quad \textrm{ as } \alpha \to 0.$$
Note that 
$$4(c_2' - c_1 ') \approx 14,$$
whereas 
$$c_2 - c_1 \approx 3.$$
We therefore clearly see that approaching the point $(\alpha, \beta) = (0,0) \in \M$ along the line $\alpha = \beta$ the gap function asymptotics are quite different from the asymptotics along trajectories of the form $\alpha = o(\beta^3)$.  

\begin{prop} \label{prop-isos}  Let $\{T_n\}$ be a sequence of triangles with diameter one and angles $\alpha_n \pi$, $\beta_n \pi$ and $\pi(1 - \alpha_n - \beta_n)$.  Assume that $\alpha_n \leq \beta_n$ and $\beta_n \to 0$ as $n \to \infty$.  Then 
$$\xi(T_n) \approx \frac{4(c_2' - c_1')}{\alpha_n^{4/3}} \quad \textrm{ as } n \to \infty \iff \alpha_n = \beta_n + o(\beta_n ^3) \textrm{ as } \beta_n \to 0.$$
\end{prop} 

\begin{proof}  By abuse of notation we shall again drop the subscript $n$.  We again estimate using domain monotonicity.  

\begin{figure}
\includegraphics[width=4in]{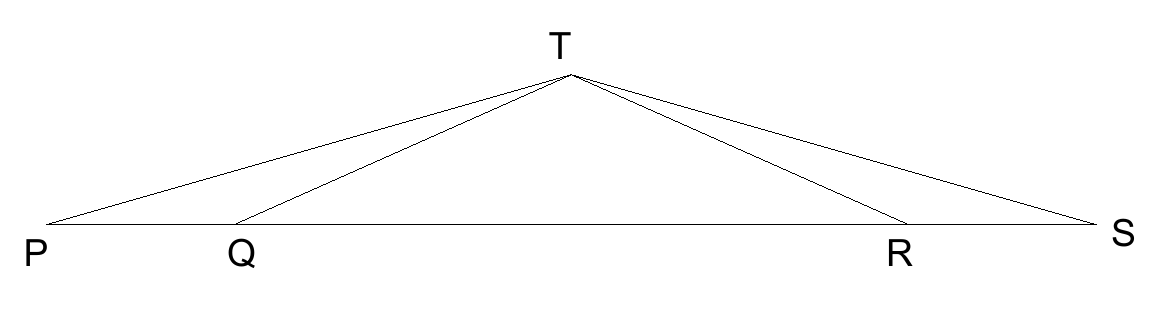}\caption{Triangle with two angles collapsing and approximating isosceles triangles.}\label{nr1}
\end{figure}

We consider the triangle $PRT$ in Figure \ref{nr1} and estimate using a smaller isosceles triangle $QRT$ and a larger isosceles triangle $PST$.  By domain monotonicity 
$$\lambda_k (PST) \leq \lambda_k (PRT) \leq \lambda_k (QRT).$$
The triangle $QRT$ has angles $\beta \pi$, $\beta \pi$, and $\pi - 2 \beta \pi$, whereas the triangle $PST$ has angles $\alpha \pi$, $\alpha \pi$, and $\pi - 2 \alpha \pi$.  By \cite{freitas}, 
$$\lambda_k (QRT) = \frac{4}{(QR)^2} \left( \frac{1}{\beta^2} + \frac{c_k '}{\beta^{4/3}} \right) + O(\beta^{-2/3}), \quad k=1,2,$$
and 
$$\lambda_k (PST) = \frac{4}{(PS)^2} \left( \frac{1}{\alpha^2} + \frac{c_k '}{\alpha^{4/3}} \right) + O(\alpha^{-2/3}), \quad k=1,2.$$
By the Law of Sines and the Law of Cosines, 
$$QR = 2A \cos(\beta \pi) = \frac{2B \sin(\alpha \pi) \cos(\beta \pi)}{\sin(\beta \pi)},$$
where $A$ and $B$ denote the sides of the triangle $PRT$ opposite the angles $\alpha \pi$, $\beta \pi$, respectively.  Since $\beta \to 0$, $B = \frac{1}{2} + O(\beta)$.  Since $\alpha \pi \leq \beta \pi$, we have 
$$QR = 1 + O(\beta).$$
Since $PQ + QR = 1$ and $PQ = RS$, it also follows that 
$$PS = 1 + O(\beta).$$
The estimate for $\xi(PRT)$ is 
$$\lambda_2 (PST) - \lambda_1 (QRT) \leq \xi(PRT) \leq \lambda_2(QRT) - \lambda_1 (PST).$$
A necessary and sufficient condition for $\xi(PRT) \sim 4(c_2' - c_1') \beta^{-4/3}$ is 
$$\frac{1}{(QR)^2 \beta^2 } - \frac{1}{(PS)^2 \alpha^2} \to 0, \quad \beta \to 0.$$
This is satisfied if and only if 
$$\alpha = \beta + o(\beta^3).$$
\end{proof}

\subsection{A general collapsing triangle}

\begin{figure}
\includegraphics[width=4in]{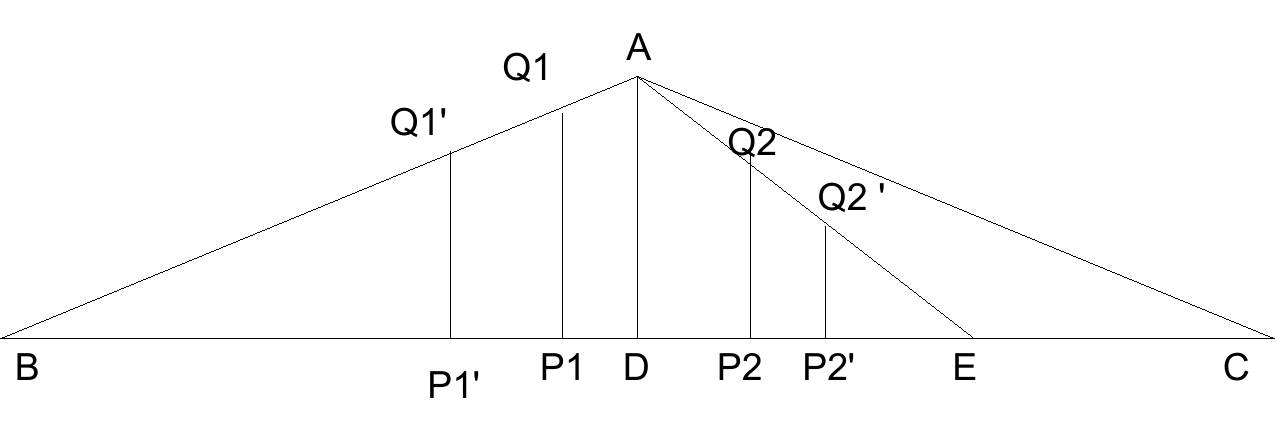}\caption{Arbitrary collapsing triangle.}\label{last}
\end{figure}

Assume the smallest angle of the triangle $\angle B = \alpha \pi;$ see Figure \ref{last}.  We wish to estimate $\lambda_2 (ABE) - \lambda_1 (ABE)$ from below.  Assume $|DE| \leq |BD| =1,$ and fix 
$$0 < \epsilon < \frac{2}{9}.$$ 
Let $P_1$ be between $B$ and $D$ so that $|P_1 D| = \alpha^{\epsilon},$ and $P_1'$ be between $B$ and $P_1$ so that $|P_1 'D| = 2 \alpha^{\epsilon}.$  If $|DE| > 2 \alpha^{\epsilon},$ let $P_2$ be between $D$ and $E$ so that $|P_2 D| = \alpha^{\epsilon},$ and $P_2'$ be between $P_2$ and $E$ so that $|P_2'D| = 2 \alpha^{\epsilon}.$  If $|DE| \leq 2 \alpha^{\epsilon},$ we do not define or use the points $P_2, P_2'.$  If $|DE| > 2 \alpha^{\epsilon}$ let $U$ be the trapezoid $A Q_1 P_1 P_2 Q_2$ and similarly let $U' = A Q_1' P_1' P_2' Q_2'.$ Let $V = ABE - U$ and let $V' = ABE - U';$ see Figure \ref{last}.  If $|DE| \leq 2 \alpha^{\epsilon},$ we let $U = A Q_1 P_1 E,$ $U' = A Q_1 ' P_1' E,$ $V = ABE - U$ and $V' = ABE - U'.$  In the estimates to follow, we show that we may estimate $\lambda_2 (ABE) - \lambda_1 (ABE)$ using $\lambda_2 (U') - \lambda_1 (U').$ 

Let $f_i$ be the eigenfunction for $\lambda_i = \lambda_i (ABE),$ $i=1,2.$  The height of $V$ is at most 
$$(1- \alpha^{\epsilon}) \tan(\alpha \pi) \approx (1-\alpha^{\epsilon}) \pi \alpha.$$
By the one dimensional Poincar\'e Inequality
$$\frac{\int_V |\nabla f_i|^2 }{\int_V f_i^2} \geq \frac{\pi^2}{(1-\alpha^{\epsilon})^2 \pi^2 \alpha^2}, \, 
\frac{\int_U |\nabla f_i|^2 }{\int_U f_i^2} \geq \frac{\pi^2}{ \pi^2 \alpha^2}, \textrm{ and } \frac{\int_{U'} |\nabla f_i|^2 }{\int_{U'} f_i^2} \geq \frac{\pi^2}{\pi^2 \alpha^2}.$$

Assume 
$$\int_{ABE} f_i ^2 = 1,  \textrm{ and } \int_V f_i ^2 = \beta_i.$$
By the variational principle, 
$$\frac{\beta}{(1- \alpha^{\epsilon})^2 \alpha^2} + \frac{1-\beta}{\alpha^2} \leq \int_V |\nabla f_i|^2 + \int_U |\nabla f_i|^2 = \lambda_i \leq \lambda_2 (ABD) = \frac{1}{\alpha^2} + \frac{c_2}{\alpha^{4/3}} + O(\alpha^{-1}).$$
Therefore,
$$\beta_i \leq \frac{\alpha^{2/3} c_2 (1- \alpha^{\epsilon})^2}{\alpha^{\epsilon} (2 - \alpha^{\epsilon})} \leq  \alpha^{2/3 - \epsilon}, \quad i=,2.$$
For simplicity in the arguments to follow, we replace all constant factors multiplying positive powers of $\alpha$ by a constant factor of $1,$ since no generality is lost as $\alpha \to 0.$   

Let $\rho$ be a smooth compactly supported function so that
\begin{equation} \label{rho1} \rho|_U \equiv 1, \quad \rho|_{V'} \equiv 0. \end{equation}
We may choose $\rho$ so that 
\begin{equation} \label{rho2} |\nabla \rho| \leq \frac{1}{\alpha^{\epsilon}}, \quad \textrm{and} \quad |\Delta \rho|, |\Delta (\rho^2)| \leq \frac{1}{\alpha^{2 \epsilon}}. \end{equation}
For the arguments to follow, we use the sign convention for the Euclidean Laplacian so that $- \Delta$ has positive spectrum.  Note that
\begin{equation} \label{last1}-(\rho f_i) \Delta (\rho f_i) = \lambda_i \rho^2 f_i ^2 - f_i ^2 \rho \Delta \rho - 2 f_i \rho (\nabla \rho )(\nabla f_i). \end{equation}

\subsubsection{Estimate for $\lambda_1 (U')$}  

Since $\rho$ vanishes on the boundary of $U',$ $\rho f_1$ is an admissible test function for the Rayleigh quotient on $U'$ (\ref{r1}) which we may use to estimate $\lambda_1 (U')$ from above.  By (\ref{last1}), 
$$\lambda_1 (U') \leq \lambda_1 + \frac{\int_{U'} \left(- \rho \Delta \rho f_1 ^2 - 2 \rho \nabla \rho f_1 \nabla f_1\right)}{\int_{U'} \rho^2 f_1 ^2}.$$
Since
$$\int_{U'} \rho \nabla \rho f_1 \nabla f_1 = \frac{1}{4} \int_{U'} \nabla \rho^2 \nabla f_1 ^2 = - \frac{1}{4} \int_{U'} \Delta \rho^2 f_1^2,$$
and $\nabla \rho, \Delta \rho = 0$ on $U,$ 
we have 
$$\int_{U'} \left(- \rho \Delta \rho f_1 ^2 - 2 \rho \nabla \rho f_1 \nabla f_1\right) \leq \frac{1}{\alpha^{2 \epsilon}} \int_{U' - U} f_1 ^2 \leq \frac{\beta}{\alpha^{2 \epsilon}} \leq \alpha^{2/3 - 3 \epsilon}.$$
Noting that
$$\int_{U'} \rho^2 f_1^2 \geq \int_U \rho^2 f_1 ^2 = 1-\beta \geq 1 - \alpha^{2/3 - \epsilon},$$
we then have
\begin{equation} \label{last-l1}\lambda_1 (U') \leq \lambda_1 + \frac{\alpha^{2/3 - 3 \epsilon}}{1 - \alpha^{2/3 - \epsilon}} \leq \lambda_1 + \alpha^{2/3 - 3 \epsilon}.  \end{equation} 
This gives the required estimate for the first eigenvalue.  

\subsubsection{Estimate for $\lambda_2(U')$}  Since $\rho f_2$ is not \`a priori orthogonal to the first eigenfunction for $U',$ we must modify it to use the Rayleigh quotient (\ref{r2}) to estimate $\lambda_2 (U').$   Since $\rho f_1$ is not orthogonal to the first eigenfunction for $U'$ because both are positive, there is some $a \in \R$ such that $\rho f_2 + a \rho f_1$ is orthogonal to the first eigenfunction for $U'.$  We may then use $\rho f_2 + a \rho f_1$ as a test function for the Rayleigh quotient on $U'.$  Integrating by parts, 
\begin{equation} \label{ip0} \int_{U'} | \nabla \rho f_i |^2 = - \int_{U'} \rho f_i \Delta(\rho f_i) = \lambda_i \int_{U'} \rho^2 f_i ^2 - \frac{1}{2} \int_{U'} \nabla \rho^2 \nabla f_i ^2 - \int_{U'} \rho \Delta \rho f_i ^2. \end{equation} 
We estimate 
$$ \left| \int_{U'} | \nabla (\rho f_i)|^2 - \lambda_i \int_{U'} \rho^2 f_i ^2 \right| \leq \frac{1}{2} \left| \int_{U'} \Delta \rho^2 f_i^2 \right| + \int_{U'} |\rho \Delta \rho | f_i^2,$$
\begin{equation} \label{last2} \leq \alpha^{-2 \epsilon} \int_{U' -U} f_i ^2 \leq \alpha^{-2\epsilon} \int_V f_i ^2 \leq \alpha^{2/3 - 3\epsilon}, \end{equation}
since $\Delta \rho$ and $\nabla \rho$ vanish identically on $U.$  We compute that $\int_{U'} \nabla(\rho f_1) \nabla(\rho f_2) = $
\begin{equation} \label{ip1} - \int_{U'} \rho f_1 \Delta(\rho f_2) = \lambda_2 \int_{U'} \rho^2 f_1 f_2 - 2 \int_{U'} \rho \nabla \rho f_1 \nabla f_2 - \int_{U'} \rho f_1 \Delta \rho f_2 \end{equation} 
and  $\int_{U'} \nabla(\rho f_1) \nabla(\rho f_2) = $
\begin{equation} \label{ip2} - \int_{U'} \rho f_2 \Delta(\rho f_1) = \lambda_1 \int_{U'} \rho^2 f_1 f_2 - 2 \int_{U'} \rho \nabla \rho f_2 \nabla f_1 - \int_{U'} \rho f_2 \Delta \rho f_1. \end{equation}  
This gives the inequality
$$ \left| 2 \int_{U'} \nabla(\rho f_1) \nabla(\rho f_2) -  (\lambda_1 + \lambda_2) \int_{U'} \rho^2 f_1 f_2 \right| \leq \int_{U'-U} \frac{|f_1 f_2| }{\alpha^{2 \epsilon}}+ 2 \left| \int_{U'-U} \rho \nabla \rho \nabla (f_1 f_2) \right| $$
\begin{equation} \label{last3} \leq \alpha^{-2 \epsilon} \left( \int_{V} |f_1 f_2| + 2 \int_{V} |f_1 f_2| \right) \leq \alpha^{2/3 - 3 \epsilon}, \end{equation}
which follows from integration by parts and the Schwarz inequality.  Expanding,
$$\int_{U'} | \nabla(\rho f_2 + a \rho f_1) |^2 - \lambda_2 \int_{U'} (\rho f_2 + a \rho f_1)^2 = I + II + III,$$
where 
$$I = \int_{U'} |\nabla \rho f_2|^2 - \lambda_2 \int_{U'} \rho^2 f_2 ^2 , \quad II =  a^2 \left( \int_{U'} |\nabla \rho f_1|^2 - \lambda_2 \int_{U'}  \rho^2 f_1 ^2 \right),$$
and 
$$ III = \int_{U'} 2a (\nabla \rho f_1)(\nabla \rho f_2) - \lambda_2 \int_{U'} 2a \rho^2 f_1 f_2.$$
By (\ref{last2}), 
\begin{equation}\label{last-1} I  \leq \alpha^{2/3 - 3 \epsilon}.  \end{equation}  
To estimate II, we note that $\lambda_1 \leq \lambda_2$ so that 
\begin{equation}\label{last-2} II \leq a^2 \left( \int_{U'} |\nabla \rho f_1|^2 - \lambda_1 \int_{U'}  \rho^2 f_1 ^2 \right) \leq a^2 \alpha^{2/3 - 3 \epsilon},\end{equation}
which follows from (\ref{last2}). 

Note that by the orthogonality of $f_1$ and $f_2$ and the Schwarz inequality, 
\begin{equation} \label{last4} \left| \int_{U'} \rho^2 f_1 f_2 \right| = \left| \int_{ABE} (1 - \rho^2) f_1 f_2 \right| \leq \left| \int_V f_1 f_2 \right| \leq \alpha^{2/3 - \epsilon}.  \end{equation}
By (\ref{last3}), (\ref{last4}), and adding and subtracting $\lambda_1 \int_{U'} a \rho^2 f_1 f_2,$ we estimate III,   
$$III \leq |a| \left| 2 \int_{U'} (\nabla \rho f_1) (\nabla \rho f_2) - (\lambda_1 + \lambda_2)  \int_{U'} \rho^2 f_1 f_2 \right| + |a| (\lambda_2 - \lambda_1) \left| \int_{U'} \rho^2 f_1 f_2 \right|$$
\begin{equation} \label{last-3} \leq |a| (\alpha^{2/3 - 3 \epsilon} + (\lambda_2 - \lambda_1) \alpha^{2/3 - \epsilon}).\end{equation}  
By (\ref{last4}), and since $\int_{U'} f_i^2 \geq 1 - \alpha^{2/3 - \epsilon},$
\begin{equation} \label{last-5} \int_{U'} (\rho f_2 + a \rho f_1)^2  \geq (1+a^2)(1-\alpha^{2/3 - \epsilon}) - 2|a| \alpha^{2/3 - 3 \epsilon} . \end{equation}
We now estimate the Rayleigh quotient for $\rho f_2 + a \rho f_1$ using (\ref{last-1}), (\ref{last-2}), (\ref{last-3}),  and (\ref{last-5})
$$\lambda_2 (U') \leq \lambda_2 + \frac{\alpha^{2/3 - 3 \epsilon} + a^2 \alpha^{2/3 - 3\epsilon} + |a|\alpha^{2/3 - 3\epsilon} + |a| (\lambda_2 - \lambda_1) \alpha^{2/3 - \epsilon}}{(1+a^2)(1-\alpha^{2/3 - \epsilon}) - 2|a| \alpha^{2/3 - 3 \epsilon}},$$
which shows that
$$\lambda_2 (U') \leq \lambda_2 + \alpha^{2/3 - 3\epsilon} + (\lambda_2 - \lambda_1)\alpha^{2/3 - \epsilon}.$$

\subsubsection{Gap estimate}

Using our estimates for $\lambda_i (U'),$ 
$$\lambda_2 - \lambda_1 \geq \lambda_2(U') - \lambda_1 (U') -\alpha^{2/3 - 3\epsilon} + (\lambda_2 - \lambda_1)\alpha^{2/3 - \epsilon} - \alpha^{2/3 - 3\epsilon}.$$
We then have 
$$(\lambda_2 - \lambda_1) \geq \frac{1}{1- \alpha^{2/3 - \epsilon} } \left(\lambda_2 (U') - \lambda_1 (U') \right) - \alpha^{2/3 - 3\epsilon},$$
which shows that 
$$\lambda_2 - \lambda_1 \geq (\lambda_2 (U') - \lambda_1 (U') ) - O(\alpha^{2/3 - 3 \epsilon}).$$
By the main theorem of \cite{swyy}, since the diameter of $U'$ is at most $4 \alpha^{\epsilon},$ 
$$\lambda_2(U') - \lambda_1 (U') \geq \frac{\pi^2}{64 \alpha^{2 \epsilon}},$$
which shows that 
$$\lambda_2  - \lambda_1 \geq C \alpha^{-2 \epsilon}$$
and is therefore unbounded as $\alpha \to 0.$  We have shown that for any triangle with one or two small angles, $\xi(T)$ becomes unbounded as the small angle(s) collapse.   \qed

\section{Collapsing polygonal domains}
Some readers may find our results for collapsing triangles counterintuitive because of the familiar example of collapsing rectangles.  For triangles, the collapse is not uniform; points opposite the longest side collapse at different rates.  Points near the vertex opposite the longest side collapse ``more slowly'' in some sense than points near the smallest angle.  For general polygonal domains, there is a subtle relationship between uniformity of collapse and the behavior of the gap.  


\begin{theorem}
\label{mgons}
Let $\{Q_n\}_{n \in \N}$ be convex $m$-gons in $\{ (x,y) \in \R^2 : x, y \geq 0 \}$.  Assume the longest side of $Q_n$, denoted $S_n = \{(x, 0) : 0 \leq x \leq 1\}$.  Let the height $h_n$ of $Q_n$ be defined by 
\begin{equation} \label{height} h_n := \inf \{ b: \exists a > 0 \textrm{ such that } Q_n \subset [0, a] \times [0, b] \}, \end{equation}
and assume $h_n \to 0$ as $n \to \infty$.  
\begin{enumerate}
\item If there exist rectangles $R_n = [0, A_n] \times [0, h_n] \supset Q_n \supset r_n = [0, a_n] \times [0, b_n]$ and a constant $\epsilon > 0$ such that $a_n \geq \epsilon$ for all $n \in \N$ and 
$$\frac{h_n}{b_n} = 1 + o(h_n^2) \quad \textrm{ as } n \to \infty,$$
then $\xi(Q_n)$ is bounded as $n \to \infty.$  
\item If there exists a convex inscribed polygon $U_n \subset Q_n$ for which the following conditions are satisfied, then $\xi(Q_n) \to \infty$ as $n \to \infty.$  
\begin{enumerate}
\item The diameter of $U_n \to 0$ as $n \to \infty.$ 
\item The longest side $\Sigma_n$ of $U_n$ is contained in $S_n.$  
\item The height of $U_n = h_n.$    
\item The height of $V_n := Q_n - U_n,$ satisfies $h_n - h(V_n) \sim h_n ^\delta$ for some $\delta \in (0, 2/9)$.
\end{enumerate}
\end{enumerate}
\end{theorem} 


\subsection{Proof of Theorem \ref{mgons}:  bounded gap}
The first case follows from domain monotonicity estimates: 
$$\lambda_2(Q_n) - \lambda_1(Q_n) \leq \lambda_2 (r_n) - \lambda_1 (R_n).$$
We may assume, since $b_n \leq h_n \to 0$ that $a_n \geq b_n$.  Then 
$$\lambda_2 (r_n) - \lambda_1 (R_n) = \pi^2 \left( \frac{h_n^2 - b_n^2}{h_n^2 b_n^2} + \frac{4}{a_n^2} - \frac{1}{A_n^2} \right).$$
Since $Q_n \subset R_n$, $A_n \geq 1$, and by assumption $a_n$ is also bounded below.  Thus 
$$\xi(Q_n) = \leq d(Q_n^2) \pi^2  \left( \frac{(h_n/b_n)^2 - 1}{h_n^2 } + \frac{4}{a_n^2} - \frac{1}{A_n^2} \right).$$
The diameter $d(Q_n) \to 1$ as $n \to \infty$ and by the assumption on $b_n$ it follows that $\xi(Q_n)$ is bounded as $n \to \infty$.  

\subsection{Proof of Theorem \ref{mgons}:  unbounded gap}

\begin{figure}[ht]
\includegraphics[width=4in]{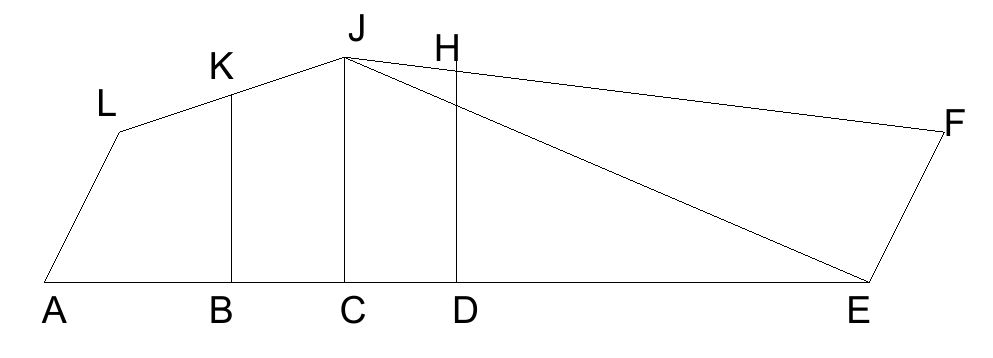}\caption{Schematic diagram of collapsing polygon with unbounded gap.}\label{gen-x}
\end{figure}

We generalize our estimates for arbitrary collapsing triangles.  Let $U' \supset U$ be a convex inscribed polygon so that one side $\Gamma$ of $U'$ satisfies $S \supset \Gamma \supset \Sigma$ and $|\Gamma - \Sigma| = h^{\delta},$  where $\delta \in (0,1)$ will be specified later.  We can define such an inscribed polygon since the diameter of $U \to 0,$ and the length of the longest side of $Q \to 1.$  Note that with these hypotheses the diameter of $U' \to 0$ as $h \to 0.$  Let $f_i$ be the eigenfunction for $Q$ with eigenvalue $\lambda_i,$ for $i=1,2.$  By convexity, $Q$ contains an inscribed right triangle $T$ of height $h = ht(U)$ and base at least $|S|/2,$ where $|S|$ is the length of the longest side of $Q.$  Scale $Q$ so that the base of $T$ is one; no generality is lost in showing $\xi$ is unbounded.  For illustration, in Figure \ref{gen-x}, $Q= AEFJL,$ $T = JCE,$ $h = |JC|,$ and $U = JKBDH.$  Assume $f_i$ are normalized so that
$$\int_Q f_i ^2 = 1, \textrm{ and let } \int_V f_i^2 = \beta_i, \quad i=1,2.$$
In the following arguments we shall absorb all constant factors independent of $h$ into a constant factor of $1$.  

By the one dimensional Poincar\'e inequality,
$$\frac{\int_V |\nabla f_i|^2 }{\int_V f_i^2} \geq \frac{\pi^2}{(ht(V))^2},$$
and 
$$\frac{\int_U |\nabla f_i|^2 }{\int_U f_i^2} \geq \frac{\pi^2}{ (ht(U))^2}, \quad \frac{\int_{U'} |\nabla f_i|^2 }{\int_{U'} f_i^2} \geq \frac{\pi^2}{ (ht(U))^2},$$
since $U$ and $U'$ have the same height.  For $h \approx 0,$ the measure of the smallest angle of $T \sim h = ht(U)$.  By domain monotonicity
$$\lambda_i (Q) \leq \lambda_2 (T) \leq \lambda_2 (S) = \frac{\pi^2}{(ht(U))^2} + \frac{c_1}{(ht(U))^{4/3}} + O(ht(U)^{-1}),$$
where $S$ is a sector of the same opening angle as $T$ which is contained in $T$ and has radius bounded below by a fixed constant.  
Estimating $\beta_i$ as we did for triangles,   
$$\frac{\beta_i \pi^2}{(ht(V))^2}+ \frac{\pi^2(1-\beta_i)}{ (ht(U))^2} \leq \int_V |\nabla f_i|^2 + \int_U |\nabla f_i|^2 = \lambda_i = \lambda_2 (T) $$
$$\leq \frac{\pi^2}{(ht(U))^2} + \frac{c_1}{(ht(U))^{4/3}} +  O(ht(U)^{-1}).$$
This gives the following estimate for $\beta_i$, 
\begin{equation} \label{beta0} \beta_i \leq \frac{ht(U)^{2/3} ht(V)^2}{ht(U)^2 - ht(V)^2}, \quad i=1,2.\end{equation}
Since by assumption $h = ht(U) - ht(V) \sim h^{1+\delta}$ for some $\delta \in (0, 2/9)$, It follows that 
\begin{equation} \label{beta1} \beta \leq h^{2/3 - \delta}. \end{equation}
Let 
\begin{equation} \label{eps} \epsilon = \frac{2}{3} - \delta. \end{equation}  

\subsubsection{Estimates for $\lambda_1$}
Define the cut-off function $\rho$ as in (\ref{rho1}) so that, 
$$|\nabla \rho | \leq h^{- \delta}, \quad | \Delta \rho | \leq h^{- 2 \delta}.$$
Using the test function $\rho f_1$ in the Rayleigh quotient for $U',$ by (\ref{ip0}),  
\begin{equation} \label{poly1} \lambda_1 (U') \leq \lambda_1(Q) + \frac{\beta}{h^{2 \delta} (1-\beta)}  \leq \lambda_1(Q) + h^{\epsilon - 2\delta}.  \end{equation}

\subsubsection{Estimates for $\lambda_2$} 
Since the function $\rho f_2$ is not \`a priori admissible as a Rayleigh quotient test function for $\lambda_2 (U')$ we again modify it to make it orthogonal to the first eigenfunction on $U'.$  So, we consider the test function 
$$\rho f_2 + a \rho f_1.$$
Using (\ref{ip0}), we estimate 
$$ \left| \int_{U'} | \nabla (\rho f_i)|^2 - \lambda_i \int_{U'} \rho^2 f_i ^2 \right| \leq \frac{1}{2} \left| \int_{U'} \Delta \rho^2 f_i^2 \right| + \int_{U'} |\rho \Delta \rho | f_i^2,$$
\begin{equation} \label{quad1} \leq h^{-2 \delta} \int_{U' -U} f_i ^2 \leq h^{-2\delta} \int_V f_i ^2 \leq h^{\epsilon- 2\delta}, \end{equation}
since $\Delta \rho$ and $\nabla \rho$ vanish identically on $U.$  By (\ref{ip1}) and (\ref{ip2}),
$$ \left| 2 \int_{U'} \nabla(\rho f_1) \nabla(\rho f_2) -  (\lambda_1 + \lambda_2) \int_{U'} \rho^2 f_1 f_2 \right| \leq \int_{U'-U} \frac{|f_1 f_2|}{h^{2 \delta}}  + 2 \left| \int_{U'-U} \rho \nabla \rho \nabla (f_1 f_2) \right| $$
\begin{equation} \label{quad2} \leq h^{-2 \delta} \left( \int_{V} |f_1 f_2| + 2 \int_{V} |f_1 f_2| \right) \leq h^{\epsilon - 2 \delta}, \end{equation}
which follows from integration by parts and the Schwarz inequality.  Expanding,
$$\int_{U'} | \nabla(\rho f_2 + a \rho f_1) |^2 - \lambda_2 \int_{U'} (\rho f_2 + a \rho f_1)^2 = I + II + III,$$
where 
$$I = \int_{U'} |\nabla \rho f_2|^2 - \lambda_2 \int_{U'} \rho^2 f_2 ^2 , \quad II =  a^2 \left( \int_{U'} |\nabla \rho f_1|^2 - \lambda_2 \int_{U'}  \rho^2 f_1 ^2 \right)$$
and 
$$ III = \int_{U'} 2a (\nabla \rho f_1)(\nabla \rho f_2) - \lambda_2 \int_{U'} 2a \rho^2 f_1 f_2.$$
By our calculations for triangles, our estimate for $\beta,$ and the estimates (\ref{quad1}) and (\ref{quad2}), 
\begin{equation} \label{q1} I \leq h^{\epsilon - 2 \delta}, \end{equation}
\begin{equation} \label{q2} II \leq a^2 h^{\epsilon - 2 \delta}, \textrm{ and } \end{equation}
\begin{equation} \label{q3} III \leq |a| (h^{\epsilon - 2 \delta } + (\lambda_2 (Q) - \lambda_1 (Q) )h^{\epsilon - 2 \delta} ).  \end{equation}
Moreover, 
\begin{equation} \label{q-lower} \int_{U'} (\rho f_2 + a \rho f_1)^2 \geq (1 + a^2) (1-h^{\epsilon}) - 2|a| h^{\epsilon - 2 \delta}. \end{equation}
Using these estimates, we estimate the Rayleigh quotient for $\rho f_2 + a \rho f_1,$
\begin{equation} \label{q-lambda2} \lambda_2 (U') \leq \lambda_2 (Q) + \frac{h^{\epsilon - 2 \delta} + a^2 h^{\epsilon - 2 \delta}+ |a| (h^{\epsilon - 2 \delta} + (\lambda_2 (Q) - \lambda_1 (Q) )h^{\epsilon - 2 \delta} )}{(1 + a^2) (1-h^{\epsilon}) - 2|a| h^{\epsilon - 2 \delta}}.  \end{equation} 
By considering the behavior as $h \to 0,$ 
$$\lambda_2 (U') \leq \lambda_2 (Q) + h^{\epsilon - 2 \delta} + h^{\epsilon - 2 \delta} (\lambda_2 (Q) - \lambda_1 (Q)).$$

\subsubsection{Gap estimate}

Using our estimates for $\lambda_i (U'),$
$$\lambda_2 (Q) - \lambda_1 (Q) \geq \lambda_2 (U') - \lambda_1 (U') - h^{\epsilon - 2 \delta}   + (\lambda_2 (Q) - \lambda_1 (Q))h^{\epsilon - 2 \delta},$$
which shows that 
\begin{equation} \label{poly-next} \lambda_2 (Q) - \lambda_1 (Q) \geq \frac{\lambda_2 (U') - \lambda_1(U')}{1 - h^{\epsilon - 2 \delta}} - h^{\epsilon - 2 \delta}.\end{equation} 
Since 
$$\epsilon = 2/3 - \delta, \quad \delta < 2/9,$$
it follows that 
$$\epsilon - 2 \delta = 2/3 - 3 \delta > 0.$$
By hypothesis on $U$ and definition of $U'$, the diameter of $U'$ vanishes as $h \to 0$, so (\ref{poly-next}) and the main theorem of \cite{swyy}  imply  
$$\lambda_2 (Q) - \lambda_1 (Q) \to \infty, \quad \textrm{ as }  h \to 0.$$ 
\qed

\begin{remark}
This technical theorem shows that the gap function is sensitive to the rate at which boundary points collapse to the base, which in two dimensions is a fixed segment.  We shall demonstrate in the following section that the same phenomenon is true for domains in $\R^n$ under one dimensional collapse.
\end{remark}

\section{The general case}
In this section, we prove the general compactness theorem:  under one dimensional collapsing, the compactness theorem is true for any convex domain.  

\begin{proof}[Proof of Theorem~\ref{main}]  In what follows, the constant $C$ is independent of $\delta$ and $\eps$ but may be different from line to line and within the same line.
Let $u_1, u_2$ be the two $\cL^2$ normalized eigenfunctions of $\Omega_\eps$ for the eigenvalues $\lambda_1(\Omega_\eps), \lambda_2(\Omega_\eps)$, respectively. Let
\begin{align*}
& E(\delta)=\{x\in E\mid h(x)>\sigma-\delta\};\\
& U(\delta)=(E(\delta)\times \mathbb R)\cap \Omega_\eps;\\
& V(\delta)=\Omega_\eps\backslash U(\delta)
\end{align*}
For $i=1,2$, by the 1-dimensional Poincar\'e inequality, we have
\begin{equation}\label{cru}
\lambda_i(\Omega_\eps)=\int_{\Omega_\eps} |\nabla u_i|^2
\geq\frac{\pi^2}{(\sigma-\delta)^2\eps^2}\int_{V(\delta)} u_i^2+\frac{\pi^2}{\sigma^2\eps^2}
\left(1-\int_{V(\delta)} u_i^2\right).
\end{equation} 
On the other hand, $\Omega_\eps$ contains a cylinder of height $(\sigma-\delta^2)$ over $E(\delta^2)$.
Thus we have
\begin{equation}\label{cru-2}
\lambda_i(\Omega_\eps)\leq \lambda_i (U(\delta^2)) \leq \lambda_i(E(\delta^2))+\frac{\pi^2}{(\sigma-\delta^2)^2\eps^2}.
\end{equation} 

Let $x_0$ be the maximum point of $h(x)$. Then by the Lipschitz continuity of $h$ there is a constant $C$, such that
\[
E(\delta^2)\supset B_{x_0}(C\delta^2).
\]
Therefore, we have
\[
\lambda_i(E(\delta^2))\leq C\delta^{-4}.
\]

Combining  the above inequality with with ~\eqref{cru},~\eqref{cru-2}, we have

\[
\frac{\pi^2}{(\sigma-\delta)^2\eps^2}\int_{V(\delta)} u_i^2+\frac{\pi^2}{\sigma^2\eps^2}
\left(1-\int_{V(\delta)} u_i^2\right)\leq \frac{C}{\delta^4}+\frac{\pi^2}{(\sigma-\delta^2)^2\eps^2},
\]
and hence  we have
\[
\int_{V(\delta)} u_i^2\leq C(\delta+\frac{\eps^2}{\delta^5}).
\]
We choose $\delta=\eps^{1/3}$. Then we have
\begin{equation}\label{45}
\int_{V(\delta)} u_i^2\leq C\eps^{1/3}.
\end{equation}

Let $\rho$ be a cut-off function which vanishes outside {$B_{x_0}(C\delta^2/2)$} with $|\nabla\rho|\le C\delta^{-2}=C\eps^{-2/3}$. 
Then
since $u_1$ is orthogonal to $u_2$, we have
\begin{equation}\label{46}
\left|\int_{U(\delta)} \rho^2 u_1u_2\right|\leq  C\eps^{1/3}.
\end{equation}

Let 
\[
\psi=u_2/u_1
\]
and let 
\[
\alpha=\frac{\int_U\psi u_1^2}{\int_U u_1^2},\quad \tilde\psi=\psi-\alpha.
\]

Let $\mu$ be the first Neumann eigenvalue of $U$ with respect to the weight $-2\log u_1$.  This is with respect to the Bakry-\'Emery Laplace operator
$$\Delta + 2\, \nabla \log u_1 \cdot \n \cdot,$$
with respect to the weighted $\cL^2$ norm $e^{2 \log u_1} dV = u_1^2 dV$, see \cite{mrl} for more details.  Then by the variational principle, we have
\[
\mu\leq \frac{\int_U|\nabla \psi|^2 u_1^2}{\int_U \tilde \psi^2\,u_1^2}
\leq \frac{\int_{\Omega_\eps} |\nabla \psi|^2 u_1^2}{\int_U \tilde \psi^2\,u_1^2} = \frac{ \lambda_2(\Omega_\eps)-\lambda_1(\Omega_\eps)}{ \int_U \tilde \psi^2\,u_1^2 }.
\]

By the above estimates for the eigenfunctions on $V$ we obtain
that for $\eps$ sufficiently small,
\[
\mu\leq 2(\lambda_2(\Omega_\eps)-\lambda_1(\Omega_\eps)).
\]

By using the same proof as in~\cite{ac}, we obtain (with the modulus of convexity of $\log u_1$ being $0$ on $U$)
\[
\mu\geq \pi^2/d(U)^2,
\]
which  completes the proof. 

\end{proof}

\begin{remark}
The assumption of the theorem, that the diameter of $U(\delta)$ goes to zero as $\delta\to 0$ is necessary. If $h(x)$ is constant, then $\Omega_\eps$ is a cylinder over $E$.  In this case, the gap converges to the gap of $E$.  More generally, if $h(x)$ is constant in a neighborhood of its maximum point, then we expect that the gap remains bounded.  
\end{remark}

The settings of Theorem \ref{main}, \cite{bf-thin}, and \cite{strip} all require the existence of a fixed, convex function $h$ (with the appropriate regularity, in our case Lipschitz continuous) such that, in two dimensions the collapsing domains are 
\begin{equation} \label{eh-collapse}\Omega_\epsilon = \{ (x,y) \mid 0 \leq x \leq 1, \quad y \leq \epsilon h(x) \}. \end{equation} 
Since $h$ is a fixed function, Theorem \ref{main} cannot be applied to the cases we considered in \S 2 and \S 3 above.  Those cases correspond to a family of domains 
$$\Omega_\epsilon = \{ (x,y) \mid 0 \leq x \leq 1, \quad y \leq h_{\epsilon} (x) \},$$
where $\{h_\epsilon\}$ is a one-parameter family of Lipschitz continuous convex functions on $[0,1]$ such that 
$$\lim_{\epsilon \to 0} h_{\epsilon} (x) \to 0, \quad \forall x \in [0,1].$$
For general collapsing geometry, there does not always exist a fixed function $h$ such that the collapse is described by a family of domains as in (\ref{eh-collapse}).  We expect that for one or higher dimensional collapse, the gap either converges to the gap of the base domain, or it diverges, and that this is determined by the uniformity with which $h_\epsilon(x) \to 0$ for different points $x$ in the base domain.  



\begin{thebibliography}{1}

\bibitem{handbook} M.  Abramowitz and I. Stegun, \em 
Handbook of mathematical functions with formulas, graphs, and
              mathematical tables, \em National Bureau of Standards Applied Mathematics, Series 55, for sale by the Superintendent of Documents, U.S. Government
              Printing Office, Washington, D.C., 1964.
               

\bibitem{ac} B. Andrews and J. Clutterbuck, \em Proof of the fundamental gap conjecture, \em J. Amer. Math Soc., 24(3), 2011, 899--916.

\bibitem{af-tri} P. Antunes and P. Freitas, \em On the inverse spectral problem for Euclidean triangles, \em Proc. R. Soc. Ser. A Math. Phys. Eng. Sci. 467, no. 2130, (2011), 1546--1562.  

\bibitem{bf-thin} D. Borisov and P. Freitas, \em Asymptotics of Dirichlet eigenvalues and eigenfunctions of the Laplacian on thin domains in $\R^d$, \em J. Funct. Anal. 258, no. 3, (2010), 893--912.  

\bibitem{chavel} I. Chavel, \em Eigenvalues in Riemannian Geometry,  \em Academic Press, (1984). 


\bibitem{cour-hil} R. Courant and D. Hilbert, \em Methods of Mathematical Physics, Volume 1, \em Interscience Publishers, (1937).  

\bibitem{durso} C. Durso, \em Solution of the inverse spectral problem for triangles, \em PhD Thesis, Massachusetts Institute of Technology, (1990).  

\bibitem{fincha} S. Finch, \em Airy Function Zeroes, \em supplementary material for \em Mathematical Constants, \em Cambridge University Press, (2003).  

\bibitem{finchb} S. Finch, \em Bessel Function Zeroes, \em supplementary material for \em Mathematical Constants, \em Cambridge University Press, (2003).  


\bibitem{freitas} P. Freitas, \em Precise bounds and asymptotics for the first Dirichlet eigenvalue of triangles and rhombi, \em J. F. A. 251, (2007), 376--398.  

\bibitem{upper} P. Freitas and D. Krejcirik, \em A sharp upper bound for the first Dirichlet eigenvalue and the growth of the isoperimetric constant of convex domains, \em Proc. Amer. Math. Soc. 136, (2008), 2997--3006.

\bibitem{strip} L. Friedlander and M. Solomyak, \em On the Spectrum of the Dirichlet Laplacian in a Narrow Strip, \em Israel J. Math, 170, (2009), 337--354.  

\bibitem{gww0} C. Gordon, D. Webb and S. Wolpert, \em One cannot hear the shape of a drum, \em Bull. Amer. Math. Soc. 27 (1992), 134--138.  

\bibitem{gww} C. Gordon, D. Webb and S. Wolpert, \em Isospectral plane domains and surfaces via Riemannian orbifolds, \em Invent. Math., \textbf{110}, no. 1, (1992), 1--22.

\bibitem{gr1} D. Grieser, \em Spectra of graph neighborhoods and scattering, \em Proc. Lond. Math. Soc. (3) 97, no. 3, (2008), 718--752. 

\bibitem{gr2} D. Grieser, \em Thin tubes in mathematical physics, global analysis and spectral geometry, \em Analysis on graphs and its applications, 565--593.  Proc. Sympos. Pure Math., 77, Amer. Math. Soc., Providence RI, (2008).  

\bibitem{grj} D. Grieser and D. Jerison, \em Asymptotics of eigenfunctions on plane domains, \em Pacific J. Math. 240 (2009), 1, 109--133.  

\bibitem{gm-phg} D. Grieser and R. Melrose, \em in preparation.  \em

\bibitem{g-tri} D. Grieser and Svenja Maronna, \em One can hear the shape of a triangle, \em preprint, arXiv:1208.3163v1 (2012).  

\bibitem{kac} M. Kac, \em Can one hear the shape of a drum? \em  Amer. Math. Monthly, \textbf{73}, no. 4, part II, (1966), 1--23.

\bibitem{ks} J. R. Kuttler and V. G. Sigillito, \em Eigenvalues of the Laplacian in two Dimensions, \em SIAM Review, 26, no. 2, (1984) 163--193.

\bibitem{l-w} T. Lang and R. Wong, \em ``Best possible'' upper bounds for the first two positive zeros of the Bessel function $J_{\nu} (x)$: the infinite case, \em J. Comp. Appl. Math., 71, (1996), 311--329. 


\bibitem{li-yau} P. Li and S. T. Yau, \em Estimates of eigenvalues of a compact Riemannian manifold. \em
Geometry of the Laplace operator (Proc. Sympos. Pure Math., Univ. Hawaii, Honolulu, Hawaii, 1979), pp. 205--239, 
Proc. Sympos. Pure Math., XXXVI, Amer. Math. Soc., Providence, R.I., (1980. 


\bibitem{l-u} L. Lorch and R. Uberti, \em ``Best possible'' upper bounds for the first two positive zeros of the Bessel functions - the finite part, \em J. Comp. Appl. Math., 75, (1996), 249--258.  

\bibitem{tri-gap} Z. Lu and J. Rowlett, \em The fundamental gap of simplices, \em to appear in Comm. in Math. Phys.  


\bibitem{mrl} Z. Lu and J. Rowlett, ] \emph{Eigenvalues of collapsing domains and drift Laplacians,} Math. Res. Lett., vol. 19, no. 3, (2012), 627--648. 


\bibitem{olver} F. W. J. Olver, \em Bessel Functions Part II, Zeros and Associated Values, \em Royal Society Math Tables Vol. 7, Cambridge University Press, (1960).  

\bibitem{polya1} G. P\'olya, \em On the characteristic frequencies of a symmetric membrane, \em Math. Z. 63, (1955), 331--337.

\bibitem{polya2} G. P\'olya and G. Szeg\"o, \em Isoperimetric inequalities in mathematical physics, \em Ann. Math. Studies, 27, Princeton Univ. Press, (1951).  

\bibitem{qw}  C. K. Qu and R. Wong, \em ``Best Possible'' Upper and Lower Bounds for the Zeros of the Bessel Function $J_{\nu} (x),$ \em Trans. of the Amer. Math. Soc., 351, no. 7, (1999), 2833--2859

\bibitem{swyy} I. M. Singer, B. Wong, S.T. Yau, and S. S. T. Yau, \em An estimate of the gap of the first two eigenvalues in the Schr\"odinger operator, \em Ann. Scuola Norm. Sup. Pisa Cl. Sci. 4, vol. 12, no. 2, (1985), 319--333.

\bibitem{bs} B. Siudeja, \em Sharp Bounds for Eigenvalues of Triangles, \em Michigan Math J., vol, 55, Issue 2, (2007), 243--254.

\bibitem{mvdb} M. van de Berg, \em On the condensation of the free Boson gas and the spectrum of the Laplacian, \em J. Stat. Phys., 31, (1983), 623--637. 

\bibitem{watson} G. N. Watson, \em A Treatise on the Theory of Bessel Functions, \em Cambridge University Press, (1966).  

\bibitem{yau1} S. T. Yau, \em An estimate of the gap of the first two eigenvalues in the
Schrodinger operator, \em Lectures on Partial Differential Equations, 223--235, 
New Stud. Adv. Math., 2, Int. Press, Somerville, MA, (2003).

\bibitem{yau} S. T. Yau, \em Nonlinear Analysis in Geometry, \em Enseignement Math., 33, (1987), 109--158.  

\bibitem{yz} Q. H. Yu and J. Q. Zhong, \em Lower bounds of the gap between the first and second eigenvalues of the Schr�dinger operator, \em Trans. Amer. Math. Soc.  294, no.1, (1986), 341--349.  

\end{thebibliography}
\end{document}